\documentclass[11pt]{article}
\usepackage[a4paper]{anysize}\marginsize{2.0cm}{2.0cm}{1.0cm}{1.0cm}
\pdfpagewidth=\paperwidth \pdfpageheight=\paperheight
\usepackage{pgf,tikz,hyperref}
\usetikzlibrary{arrows}
\usepackage{amsfonts,amssymb,amsthm,amsmath,eucal}
\usepackage{graphicx}
\usepackage{ltablex}
\usepackage{caption, booktabs,verbatim}

\theoremstyle{plain}
\newtheorem{thm}{Theorem}[section]

\newtheorem{lemma}[thm]{Lemma}

\newtheorem{corollary}[thm]{Corollary}

\theoremstyle{definition}
\newtheorem{definition}[thm]{Definition}
\newtheorem{remark}[thm]{Remark}

\newtheorem{problem}[thm]{Problem}

\newtheorem{thevarthm}[thm]{\varthmname}

\newenvironment{varthm*}[1]{\trivlist\item[]{\bf #1.}\it}{\endtrivlist}

\def\cftil#1{\ifmmode\setbox7\hbox{\accent"5E#1}\else\setbox7\hbox{\accent"5E#1}\penalty 10000\relax\fi\raise 1\ht7\hbox{\lower1.1ex\hbox to 
1\wd7{\hss\accent"7E\hss}}\penalty 10000\hskip-1\wd7\penalty 10000\box7 }

\newcommand\eps{\varepsilon}
\newcommand\be{\begin{eqnarray*}}
\newcommand\ee{\end{eqnarray*}}

\newcommand\C{\mathbb C}

\newcommand\Z{\mathbb Z}

\renewcommand\P{\mathbb P}

\newcommand\newop[2]{\def#1{\mathop{\rm #2}\nolimits}}
\newop\edim{edim}
\newop\Zeroes{Zeroes}
\newop\Jac{Jac}
\newop\Ass{Ass}
\newop\SL{SL}
\newop\PGL{{\P}GL}
\newop\Km{Km}
\newop\reg{reg}
\newop\Hess{Hess}
\newop\tr{tr}

\newcommand\keywords[1]{{\renewcommand\thefootnote{}\footnotetext{\textit{Keywords:} #1.}}}
\newcommand\subclass[1]{{\renewcommand\thefootnote{}\footnotetext{\textit{Mathematics Subject Classification (2010):} #1.}}}

\begin{document}

\author{Tomasz Szemberg, Justyna Szpond}
\title{Sextactic points on the Fermat cubic curve and arrangements of conics}
\date{\today}
\maketitle
\thispagestyle{empty}

\begin{abstract}
   The purpose of this note is to report, in narrative rather than rigorous style, about the nice geometry of $6$-division points on the Fermat cubic 
   $F$ and various conics naturally attached to them. Most facts presented here were derived by symbolic algebra programs and the idea of the note is 
   to propose a research direction for searching for conceptual proofs of facts stated here and their generalisations. Extensions in several directions 
   seem possible (taking curves of higher degree and contact to $F$, studying higher degree curves passing through higher order division points on $F$, 
   studying curves passing through intersection points of already constructed curves, taking the duals etc.) and we hope some younger colleagues might 
   find pleasure in following proposed paths as well as finding their own.
\end{abstract}

\keywords{osculating curves, division points, Fermat cubic}
\subclass{MSC 14C20 \and MSC 14N20 \and MSC 13A15}


\section{Introduction}
   Any complex elliptic curve $E$ can be embedded as a smooth cubic in the projective plane $\P^2$. With the embedding some points on $E$ are 
   distinguished by the way they interact with certain other curves in the plane. For example, it is well-known that there are exactly $9$ flex points 
   on $E$, i.e., points where the tangent line is hyper-tangent by what we mean that it intersects $E$ with multiplicity $3$ rather than the 
   multiplicity $2$ as an ordinary tangent line does. This was observed already in the 18th century by Maclaurin.

   From the intrinsic point of view choosing one of the inflexion points as the origin (in the group law on $E$) the set of all $9$ flex points is 
   exactly the set of $3$-division points $E[3]$ on $E$, i.e., points $P$ subject to the condition $3P=0$ in the addition law on $E$. This link between 
   the intrinsic and extrinsic geometry is justified by Abel's Theorem \cite[Theorem IV.4.13B]{Har77}.

   Since flex points of a smooth curve are well known to be common zeroes of the equation $f$ defining this curve and the Hessian $H(f)=H_1(f)$, this 
   provides an extrinsic, concrete way to obtain the ideal of all $3$-division points $E[3]$ on $E$ as a complete intersection ideal generated by $f$ 
   and $H(f)$.

   In the middle of the 19th century Cayley \cite{Cay59}, \cite{Cay65} observed that the interplay between hyperosculating lines (flex lines) and 
   $3$-division points can be extended in a natural manner. More specifically, points admitting hyperosculating curves of higher degree are precisely 
   higher order division points on the cubic curve. In the present note we focus on hyperosculating conics. Our interest in these curves and their 
   arrangement is partly motivated by the considerable grow of interest in arrangements of curves of higher degree in the last couple of years, see 
   e.g. \cite{PRS17}, \cite{PS20}, \cite{DLPU20}, \cite{DJP22}, \cite{KRS22}. This interest in turn has roots in various deep problems in 
   combinatorics, singularities theory and construction of surfaces of general type with certain prescribed invariants.

   Throughout the note, by a slight abuse of notation, we denote by the same symbol, say $\Gamma$ a homogeneous polynomial in $\C[x,y,z]$ and its set 
   of zeroes, i.e., the plane curve it defines. We hope that it will be clear from the context if we have an algebraic or geometric object in mind.
\section{Flex points and flex lines}
   We recall briefly properties of flex lines and determine them explicitly in the case of the Fermat cubic
   $$F=x^3+y^3+z^3.$$
   Its Hessian is readily computed (up to a scalar) as
   $$H(F)=xyz,$$
   i.e., the set of zeroes of the Hessian splits in three coordinate lines. Every flex point lies on one of these lines. They are evenly distributed: 
   three flex points on each of the coordinate lines.
   Explicit coordinates of these points are:
$$\begin{array}{ccccccccccc}
     P_1 & = & (1:-1:0), && P_2 & = & (1:-\eps:0), && P_3 & = & (1:-\eps^2:0),\\
     P_4 & = & (1:0:-1), && P_5 & = & (1:0:-\eps), && P_6 & = & (1:0:-\eps^2),\\
     P_7 & = & (0:1:-1), && P_8 & = & (0:1:-\eps), && P_9 & = & (0:1:-\eps^2),
\end{array}$$
   where $\eps$ is a primitive root of $1$ of order $3$.

   Since the lines tangent to $F$ at flex points have the order of tangency $3$, any line passing through two of flex points necessarily meets $F$ in a 
   third flex (it must be a flex by Abel's Theorem) different from the two flexes determining the line. Thus there are altogether $12$ lines meeting 
   $F$ in three mutually distinct flex points. This gives rise to a
   $$(12_3, 9_4)$$
   arrangement, which is known as the Hesse arrangement. It is an arrangement enjoying many properties interesting in various branches of mathematics.
   We just mention here a recent work of Bassa and \"Ozg\"ur K\.i\c{s}\.isel \cite{Bas20}, where they claim the only complex $4$-net is the Hesse. The 
   equations of the $12$ lines can be given explicitly:
$$\begin{array}{ccccccccccc}
     L_1 & = & x, && L_2 & = & y, && L_3 & = & z,\\
     L_4 & = & x+y+z, && L_5 & = & x+\eps y+z, && L_6 & = & x+\eps^2 y+z,\\
     L_7 & = & x+y+\eps z, && L_8 & = & x+y+\eps^2 z, && L_9 & = & x+\eps y+\eps z\\
     L_{10} & = & x+\eps y+\eps^2 z, && L_{11} & = & x+\eps^2 y+\eps z, && L_{12} & = & x+\eps^2 y+\eps^2 z,
\end{array}$$
   Dualizing, we obtain a
   $$(9_4, 12_3)$$
   arrangement of $9$ lines with a peculiar property that whenever two of them intersect there is a third one passing through the intersection point.
   Thus there are only triple intersection points among the lines forming this configuration. It is one of very few known configurations of this type, 
   see \cite[Open Problem 1.16]{Har18}. It is an interesting and challenging problem to either show that the list in \cite[Remark 1.1.4]{Har18} is 
   complete or to construct new examples.

   The $12$ points in which the lines of the dual Hesse arrangement intersect form a configuration very interesting from the algebraic point of view.
   Let $I$ be the saturated ideal of these $12$ points. Then its third symbolic power $I^{(3)}$ fails to be contained in the second ordinary power 
   $I^2$. The witness for the non-containment is the product of equations of the $9$ configuration lines, i.e., the polynomial
$$(x^3-y^3)(y^3-z^3)(z^3-x^3).$$
This property of $I$ was discovered in 2013 by Dumnicki, Tutaj-Gasi\'nska and the first author \cite{DST13}. This finding led to considerable 
development of the containment theory between symbolic and ordinary powers of homogeneous ideals reflected in a large and still growing number of 
articles on this subject \cite{Ake17, BGHN22, BCH14, real, DDGHN18, Dum15, DumTut17, HarSec15, MalSzp17ERA, MalSzp18, Sec15}.
   See our previous joint work \cite{SzeSzp17} for an introduction to this circle of ideas.
\section{Sextactic points and derived objects}
   Cayley proved in \cite{Cay59} that for a smooth point on a plane curve of degree $\geq 3$ there exists a unique osculating conic, whose local 
   intersection multiplicity with the curve at the point of contact is at least $5$.
   This is a degree $2$ analogy of a tangent line, which is distinguished among all lines passing through a smooth point of a plane curve by having the 
   local intersection multiplicity with the curve at least $2$.
   As we saw in the previous section, in some points of the curve the tangent lines have an even greater contact, namely the local intersection 
   multiplicity is $\geq 3$. Points with this property are of course the flex points of the curve and they can be identified for a smooth plane curve 
   as the intersection points of the curve with its Hessian. Thus an irreducible curve of degree $d$ has at most $3d(d-2)$ flexes.

   In \cite{Cay65} Cayley realised that a plane curve, in analogy to flexes, possesses points where the osculating conic has contact $\geq 6$.
   He introduced the terminology of \emph{sextactic points} and named accordingly the hypo-osculating conics as \emph{sextactic conics}.
   Note that at the flexes of the curve, the sextactic conic is easily identified as the tangent line taken twice.

   It is natural to wonder if the set of sextactic points of a plane curve $\Gamma$ can be identified by intersecting $\Gamma$ with some other curve 
   derived out of its equation. In order to answer this question
   Cayley introduced the notion of the second Hessian of a curve and proved that the sextactic points are, in analogy to the flexes, the common points 
   of the curve and its second Hessian, which defines a curve of degree $12\deg(\Gamma)-27$.
   Only recently Maugesten and Moe \cite{MM19} checked Cayley's formula carefully and pinned down an inaccuracy in its coefficients. Since the 
   erroneous formula was repeated in the literature for over 150 years, we state it here in the correct form following the notation of \cite{MM19}. The 
   inaccuracy of Cayley was to write 40 in the place of the coefficient 20 appearing in the third line of the formula in Definition \ref{def:H2}. We 
   refer also to our Singular script \cite{Sing1} for a verification of the formula.
\begin{definition}[The second Hessian, Cayley 1865, Maugesten, Moe 2019]\label{def:H2}
   Let $\Gamma$ be a plane curve of degree $d$. Then, its second Hessian $H_2(\Gamma)$ is
\begin{align*}
    H_2(\Gamma) & =(12d^2-54d+57)H(\Gamma)\Jac{(\Gamma, H(\Gamma), \Omega_H)}\\
    & + (d-2)(12d-27)H(\Gamma)\Jac{(\Gamma, H(\Gamma), \Omega_\Gamma)}\\
    & - 20(d-2)^2 \Jac{(\Gamma, H(\Gamma),\Psi)}.
\end{align*}
\end{definition}
   Of course now the explanation of symbols  appearing in Definition \ref{def:H2} is due. We proceed step by step.
   For polynomials $f,g,h$ we have
   $$\Jac{(f,g,h)}=\det
   \begin{pmatrix}
      f_x & f_y & f_z \\
      g_x & g_y & g_z \\
      h_x & h_y & h_z \\
   \end{pmatrix},$$
   where $f_u$ denotes as usual the partial derivative of $f$ with respect to $u$. The Hesse matrix of a polynomial $\Gamma$ is thus
$$\Hess(\Gamma)=
\begin{pmatrix}
   \Gamma_{xx} & \Gamma_{xy} & \Gamma_{xz} \\
   \Gamma_{yx} & \Gamma_{yy} & \Gamma_{yz} \\
   \Gamma_{zx} & \Gamma_{zy} & \Gamma_{zz} \\
\end{pmatrix}.$$
Let us abbreviate its determinant by $H=H(\Gamma)=\det(\Hess(\Gamma))$.
Then $H$ is itself a polynomial of degree $3(\deg(\Gamma)-2)$ and we can consider its Hesse matrix
$$\Hess(H)=
\begin{pmatrix}
   H_{xx} & H_{xy} & H_{xz} \\
   H_{yx} & H_{yy} & H_{yz} \\
   H_{zx} & H_{zy} & H_{zz} \\
\end{pmatrix}.$$
Then
$$\Omega=\tr\left(\Hess(\Gamma)^{adj}\cdot\Hess(H)\right)$$
or more explicitly $\Omega$ is the scalar product
\begin{align*}
\Omega=&
\begin{pmatrix}
   \Gamma_{yy}\Gamma_{zz}-\Gamma_{yz}^2\\
   \Gamma_{xx}\Gamma_{zz}-\Gamma_{xz}^2\\
   \Gamma_{xx}\Gamma_{yy}-\Gamma_{xy}^2\\
   \Gamma_{xy}\Gamma_{xz}-\Gamma_{xx}\Gamma_{yz}\\
   \Gamma_{xy}\Gamma_{yz}-\Gamma_{yy}\Gamma_{xz}\\
   \Gamma_{xz}\Gamma_{yz}-\Gamma_{zz}\Gamma_{xy}
\end{pmatrix}^T\cdot
\begin{pmatrix}
   H_{xx}\\
   H_{yy}\\
   H_{zz}\\
   2H_{yz}\\
   2H_{xz}\\
   2H_{xy}
\end{pmatrix}.
\end{align*}
For a variable $u\in\left\{x,y,z\right\}$ we have then finally
\begin{align*}
(\Omega_\Gamma)_u=&
\begin{pmatrix}
   (\Gamma_{yy}\Gamma_{zz}-\Gamma_{yz}^2)_u\\
   (\Gamma_{xx}\Gamma_{zz}-\Gamma_{xz}^2)_u\\
   (\Gamma_{xx}\Gamma_{yy}-\Gamma_{xy}^2)_u\\
   (\Gamma_{xy}\Gamma_{xz}-\Gamma_{xx}\Gamma_{yz})_u\\
   (\Gamma_{xy}\Gamma_{yz}-\Gamma_{yy}\Gamma_{xz})_u\\
   (\Gamma_{xz}\Gamma_{yz}-\Gamma_{zz}\Gamma_{xy})_u
\end{pmatrix}^T\cdot
\begin{pmatrix}
   H_{xx}\\
   H_{yy}\\
   H_{zz}\\
   2H_{yz}\\
   2H_{xz}\\
   2H_{xy}
\end{pmatrix}.
\end{align*}
and
\begin{align*}
(\Omega_H)_u=&
\begin{pmatrix}
   \Gamma_{yy}\Gamma_{zz}-\Gamma_{yz}^2\\
   \Gamma_{xx}\Gamma_{zz}-\Gamma_{xz}^2\\
   \Gamma_{xx}\Gamma_{yy}-\Gamma_{xy}^2\\
   \Gamma_{xy}\Gamma_{xz}-\Gamma_{xx}\Gamma_{yz}\\
   \Gamma_{xy}\Gamma_{yz}-\Gamma_{yy}\Gamma_{xz}\\
   \Gamma_{xz}\Gamma_{yz}-\Gamma_{zz}\Gamma_{xy}
\end{pmatrix}^T\cdot
\begin{pmatrix}
   (H_{xx})_u\\
   (H_{yy})_u\\
   (H_{zz})_u\\
   (2H_{yz})_u\\
   (2H_{xz})_u\\
   (2H_{xy})_u
\end{pmatrix}.
\end{align*}
Finally $\Psi$ is defined as
$$\Psi=-\det
\begin{pmatrix}
   0 & H_x & H_y & H_z\\
   H_x & \Gamma_{xx} & \Gamma_{xy} & \Gamma_{xz}\\
   H_y & \Gamma_{xy} & \Gamma_{yy} & \Gamma_{yz}\\
   H_z & \Gamma_{xz} & \Gamma_{yz} & \Gamma_{zz}
\end{pmatrix}.$$

\section{Sextactic points and conics}
   Applying Definition \ref{def:H2} to the Fermat cubic $F$ (we omit the dull calculations for the check of which we refer to our Singular script 
   \cite{Sing1}) we obtain up to a scalar
   $$H_2(F)=(x^3-y^3)(y^3-z^3)(z^3-x^3).$$
   Thus, similarly as in the case of the ordinary Hessian, the second Hessian splits completely into linear factors. The lines defined by these linear 
   forms are arranged in the Fermat configuration of order $3$ or equivalently in the dual Hesse arrangement, which we have already encountered in the 
   previous section! Together with the Hessian, we have
   $$H_1(F)\cdot H_2(F)=xyz(x^3-y^3)(y^3-z^3)(z^3-x^3),$$
   which is the \emph{extended Fermat arrangement}. For more details on exciting properties of Fermat arrangements we refer to \cite{Szp19c}.

   The set of all sextactic points on the Fermat cubic is thus a complete intersection of the cubic $F$ and its second Hessian $H_2(F)$. They can be 
   individually computed explicitly and we obtain the following list:

\renewcommand*{\arraystretch}{1.5}
$$
\begin{array}{lll}
S_1=\left[-\frac{1}{2} \varepsilon \mu^2:-\frac{1}{2} \varepsilon \mu^2:1\right], &
S_2=[1: (\varepsilon + 1) \mu : 1], &
S_3=[(\varepsilon + 1) \mu : 1 : 1],\\

S_4=\left[\frac{1}{2} (\varepsilon + 1) \mu^2 : -(\frac{1}{2}) \mu^2 : 1\right], &
S_5=[\varepsilon : -\varepsilon \mu : 1], &
S_6=[-\mu : -\varepsilon - 1 : 1],\\

S_7=\left[-\frac{1}{2} \mu^2 : \frac{1}{2} (\varepsilon + 1) \mu^2 : 1\right], &
S_8=[-\varepsilon - 1 : -\mu : 1], &
S_9=[-\varepsilon \mu : \varepsilon : 1],\\

S_{10}=\left[-\frac{1}{2} \mu^2 : -\frac{1}{2} \mu^2 : 1\right], &
S_{11}=[1 : -\mu : 1], &
S_{12}=[-\mu : 1 : 1],\\

S_{13}=\left[\frac{1}{2} \varepsilon \mu^2 : \frac{1}{2} (\varepsilon + 1) \mu^2 : 1\right], &
S_{14}=[\varepsilon : (\varepsilon + 1) \mu : 1], &
S_{15}=[-\varepsilon \mu : -\varepsilon - 1 : 1],\\

S_{16}=\left[\frac{1}{2} (\varepsilon + 1) \mu^2: -\frac{1}{2} \varepsilon \mu^2 : 1\right], &
S_{17}=[-\varepsilon - 1 : -\varepsilon \mu : 1], &
S_{18}=[(\varepsilon + 1) \mu : \varepsilon : 1],\\

S_{19}=\left[\frac{1}{2} (\varepsilon + 1) \mu^2 : \frac{1}{2} (\varepsilon + 1) \mu^2 : 1\right], &
S_{20}=[1 : -\varepsilon \mu : 1], &
S_{21}=[-\varepsilon \mu : 1 : 1],\\

S_{22}=\left[-\frac{1}{2} \mu^2 : -(\frac{1}{2}) \varepsilon \mu^2 : 1\right], &
S_{23}=[\varepsilon : -\mu : 1], &
S_{24}=[(\varepsilon + 1) \mu : -\varepsilon - 1 : 1],\\

S_{25}=\left[-\frac{1}{2} \varepsilon \mu^2 : -(\frac{1}{2}) \mu^2 : 1\right], &
S_{26}=[-\varepsilon - 1 : (\varepsilon + 1) \mu : 1], &
S_{27}=[-\mu : \varepsilon : 1],
\end{array}$$
where $\varepsilon$ is as above a primitive root of $1$ and $\mu$ is the real third root of $2$.
\begin{remark}
   It is worth to point out that the points $S_1,\ldots,S_{27}$ listed above are nothing other but the $6$-division points on $F$, which are not 
   $3$-division points.
\end{remark}
Endowing $F$ with a level $6$ structure, i.e., fixing a group isomorphism
$$\alpha: F[6]\to \Z_6\times\Z_6$$
the $3$-division points $F[3]$ correspond to pairs of even integers $(a,b)$ with $0\leq a,b\leq 4$.

Since the flex points are aligned in an unusual way, it is natural to expect that  the sextactic points behave in an unexpected way with respect to 
conics. Given explicit coordinates of points $S_1,\ldots,S_{27}$ and running a brute-force computer calculation we are able to enumerate conics passing 
through six-tuples of the sextactic points.
Note that $6$ is the maximal number of sextactic points which can be contained in a conic. Indeed, any conic through $7$ or more of sextactic points 
would be a component of $F$ by Bezout's Theorem. Of course the Fermat cubic is smooth and so, in particular, irreducible so that this situation cannot 
happen.
\begin{lemma}\label{lem:conics_through_6}
   There are exactly $8\;244$ conics containing six of the sextactic points.
\end{lemma}
\begin{proof}
   Direct computer computation, see \cite{Sing1}.
   Alternatively one can use the isomorphism $\alpha$ mentioned above. Then the task amounts to identifying $6$-tuples of pairs of integers (with at 
   least one of them odd)
   $$(a_1,b_1),\ldots,(a_6,b_6)\in \Z_6\times \Z_6$$
   with the property
   $$\sum_{i=1}^6 a_i=0 \mod(3) \;\mbox{ and }\;
   \sum_{i=1}^6 b_i=0 \mod(3).$$
   Indeed, it follows from the Abel Theorem that $6$ points on an elliptic curve are contained in conic if and only if their sum (in the addition law 
   on the elliptic curve) is a $2$-division point.
\end{proof}
\begin{corollary}\label{cor:conics_through_point}
   Since the configuration of the sextactic points is symmetric, it follows that there are exactly $1\;832$ conics among those from Lemma 
   \ref{lem:conics_through_6} passing through a fixed sextactic point.
\end{corollary}
\begin{remark}
Unlike in the case of flex lines, which passing through two flex points must pass through a third one, there are many conics passing just through $5$ 
of the sextactic points. However, all such conics do not have any additional intersection points with the Fermat cubic. It turns out that all such 
conics are tangent at one of the sextactic points to $F$.
\end{remark}
Among the $8\;244$ conics mentioned above, there is a number of reducible and a number of irreducible ones. More precisely the following holds.
\begin{lemma}\label{lem:smooth_conics_through_6}
   There are exactly $5\; 976$ \emph{smooth} conics containing six of the sextactic points.
\end{lemma}
\begin{proof}
   The proof is computational. Alternatively one can identify those conics which split into lines. To this end one can first identify the $81$ lines 
   from Lemma \ref{lem:arrangement_lines_points}. The argument here is in turn a  simplified version of the argument presented in the proof of Lemma 
   \ref{lem:conics_through_6}. Collinear triples of sextactic points correspond under the isomorphism $\alpha$ to pairs $(a_1,b_1), (a_2, b_2)\in 
   \Z_6\times \Z_6$ such that at least one of the sums $a_1+a_2$ or $b_1+b_2$ is odd.
\end{proof}
In the analogy to Corollary \ref{cor:conics_through_point} we have now.
\begin{corollary}\label{cor:smooth_conics_through_point}
   There are $1\;328$ smooth conics among those from Lemma \ref{lem:smooth_conics_through_6} passing through a fixed sextactic point.
\end{corollary}
Putting Lemma \ref{cor:conics_through_point} and Lemma \ref{cor:smooth_conics_through_point} together we obtain the following result.
\begin{corollary}\label{cor:conics_splitting_into_lines}
   There are $2\;268$ conics among those from Lemma \ref{lem:conics_through_6} splitting into two lines. Moreover, there are exactly $3$ sextactic 
   points on each of the lines and the intersection points of the lines forming such a conic are not in the set of the sextactic points.
\end{corollary}
With a little more effort we arrive to a rather surprising fact that there are altogether only $81$ lines of which reducible conics are composed. More 
precisely we obtain the following arrangement.
\begin{lemma}\label{lem:arrangement_lines_points}
   The $81$ lines and the $27$ sextactic points form a $(81_3,27_9)$ arrangement, i.e., each line contains $3$ points and through each point there are 
   $9$ lines passing.
\end{lemma}
As a byproduct from Lemma \ref{lem:arrangement_lines_points} we derive the following somewhat surprising fact.
\begin{lemma}
   The product of the equations of all $81$ lines from Lemma \ref{lem:arrangement_lines_points} is defined over $\Z$.
\end{lemma}
\section{Perspectives}
   It has been observed by Gattazzo \cite{Gat79} that a smooth plane cubic contains points which are special from the perspective of unexpectedly high 
   contact with curves of higher degree. More precisely he defines points of type $3k$ for $k\geq 1$ as follows.
\begin{definition}[Type $3k$ points]
   Let $\Gamma$ be a smooth plane cubic. A point $P\in\Gamma$ is of type $3k$, if there exists an \emph{irreducible} curve of degree $k$ intersecting 
   $\Gamma$ in $P$ with multiplicity $3k$.
\end{definition}
   Thus flex points are type $3$ points and the sextactic points are type $6$ points.

   Gattazzo observed that there are $72$ points of type $9$ on a smooth plane cubic. We finish this note with the following challenges.
\begin{problem}
   Check if there exists a curve of degree $24$ which intersects a smooth elliptic curve in type $9$ points.
\end{problem}
\begin{problem}
   Compute explicit coordinates of type $9$ points on the Fermat cubic.
\end{problem}

\paragraph*{Acknowledgement.}
Our research was partially supported by National Science Centre, Poland, Opus Grant 2019/35/B/ST1/00723.








\bigskip
\bigskip
\small

\bigskip
\noindent
   Tomasz Szemberg,\\
   Department of Mathematics, Pedagogical University of Cracow,
   Podchor\c a\.zych 2,
   PL-30-084 Krak\'ow, Poland.

\nopagebreak
\noindent
   \textit{E-mail address:} \texttt{tomasz.szemberg@gmail.com}\\

\bigskip
\noindent
   Justyna Szpond,\\
   Institute of Mathematics,
   Polish Academy of Sciences,
   \'Sniadeckich 8,
   PL-00-656 Warszawa, Poland

\nopagebreak
\noindent
   \textit{E-mail address:} \texttt{szpond@gmail.com}\\


\end{document}